\documentclass[12pt]{amsart}

\def\checkbox{\leavevmode\vbox to 9pt{\hrule \vss
	\hbox to 9pt{\vrule height 9pt \hfil\vrule height 9pt}\vss
	\hrule}\ }

\newcommand{\Q}{\mathbb Q}
\newcommand{\R}{\mathbb R}

\renewcommand{\epsilon}{\varepsilon}
\renewcommand{\phi}{\varphi}

\newtheorem{Lemma}{Lemma}
\newtheorem{Theorem}[Lemma]{Theorem}

\begin{document}

\address{Department of Mathematics, Pennsylvania State University, University Park, PA $16802$, USA.}
\author{Alexander Borisov}
\title{On a Question of Craven and a Theorem of Belyi}

\begin{abstract}
In this elementary note we prove that a polynomial with rational coefficients divides the derivative of some polynomial which splits in $\Q$ if and only if all of its irrational roots are real and simple. This provides an answer to a question posed by Thomas Craven. Similar ideas also lead to a variation of the proof of Belyi's theorem that every algebraic curve defined over an algebraic number field admits a map to $P^1$ which is only ramified above three points. As it turned out, this variation was noticed previously by G. Belyi himself and Leonardo Zapponi.
\end{abstract}

\subjclass[2000]{Primary 11R80, secondary 11G99}
\email{borisov@math.psu.edu}

\maketitle

The first goal of this short note is to answer a question posed by Thomas Craven in \cite{Craven}. To explain the question we need to recall some definitions. 
A field $K$ is called split Rolle (or for brevity SR) field if for any polynomial $f(x)\in K[x]$ which splits completely in $K$ its derivative $f'(x)$ also splits completely in $K$. For any field $L$ there is the unique (up to $L-$isomorphism) smallest SR-field $K$ containing
$L$. It can be obtained by successively adjoining the roots of derivatives of totally split polynomials. We will call this field the SR-closure of $L$.
The term ``split Rolle" is due to Craven, but the problem of classifying SR fields was actually posed by Kaplansky (cf. \cite{Kaplansky}, p.30).

In \cite{Craven} Thomas Craven determined the SR-closure for all finite fields. It turned out that except for $F_2$ and $F_4$ the SR-closure coincides with the algebraic closure. A very simple proof of this was recently found by C. Ballantine and J. Roberts in \cite{Monthly}. At the end of \cite{Craven}, T. Craven noted that it is unknown what the SR-closure of $\Q$ is. It is easy to show that it is contained in the totally real algebraic closure of $\Q,$ which is the field consisting of roots of polynomials with rational coefficients that are split in $\R$.  He wrote that one should not expect these two fields to be the same. The next theorem shows that they do coincide. Moreover, one iteration of the adjoining process is enough.

\begin{Theorem}
Suppose $f(x)\in {\Q}[x]$ is an irreducible polynomial. Then the following two conditions are equivalent:
\begin{enumerate}
\item{} All roots of $f(x)$ are real;
\item{} There exists a polynomial $F(x) \in {\Q}[x]$ such that all roots of $F(x)$ are rational and $f(x)$ divides its derivative $F'(x)$.
\end{enumerate}
\end{Theorem}

{\bf Proof.} (2)$\Rightarrow $(1). By Rolle's theorem $F'$ has a root between any two consecutive roots of $F.$ For each root of $F$ its order in $F'$ is one less than its order in $F$. The degree of $F'$ is one less than the degree of $F.$ So all the roots required by Rolle's theorem are simple, and together with the multiple roots of $F$ they exhaust all the roots of $F'.$ Since $f$ divides $F'$ all roots of $f$ are real.

(1)$\Rightarrow $(2). Suppose the degree of $f$ is
$n$, and the leading coefficient $a_n$ of $f$ is positive. Suppose
$$\alpha_1< \alpha_2< \dots <\alpha_n$$
are the roots of $f.$ Choose rational numbers $\{q_i\}$ so that
$$q_1<\alpha_1< q_2 <\alpha_2<q_3<a_3< ...<q_n<\alpha_n< q_{n+1}.$$

Then
$$f(x)=\sum \limits_{i=1}^{n+1}f(q_j)\frac{\prod_{j\neq i}(x-q_j)}{\prod_{j\neq i}(q_i-q_j)}$$
by the Lagrange interpolation formula. By construction all numbers
$$\frac{f(q_i)}{\prod_{j\neq i}(q_i-q_j)}=a_n\frac{\prod_{j=1}^{n}(q_i-\alpha_j)}{\prod_{j\neq i}(q_i-q_j)}$$
are positive.
Multiplying by a common denominator, we can express some multiple of $f(x)$  as  
$$\sum \limits_{i=1}^{n+1} k_i \prod \limits_{j\neq i}(x-q_j)$$
for some {\it positive} integers $k_1, k_2, ..., k_{n+1}.$ Now to complete the proof we notice that the above polynomial divides the
derivative of the polynomial
$$F(x)=\prod_{i=1}^{n+1}(x-q_i)^{k_i}.$$

The above argument can be generalized to completely describe the polynomials that divide a derivative of a polynomial with rational roots. 
\begin{Theorem} Suppose $f(x)\in {\Q}[x]$ is any polynomial. Then the following two conditions are equivalent:
\begin{enumerate}
\item All irrational roots of $f$ are real and simple;
\item There exists a polynomial $F(x) \in {\Q}[x]$ such that all roots of $F(x)$ are rational and $f(x)$ divides its derivative $F'(x)$.
\end{enumerate} 
\end{Theorem}

{\bf Proof.} (2)$\Rightarrow $(1). This follows from the argument at the beginning of the proof of Theorem 1.

(1)$\Rightarrow $(2). Multiplying $f$ by some polynomial with rational roots, we can assume that set of the roots of $f$ is a disjoint union of $\{q_i\}$ and $\{\alpha_i\},$ where the $q_i$ are rational and the $\alpha_i$ are simple, and 
$$q_1<\alpha_1< q_2 <\alpha_2<q_3<a_3< ...<q_n<\alpha_n< q_{n+1}.$$
Then we proceed as in the proof of the Theorem 1 to get a polynomial $F(x)$ such that the roots of $F$ are all $q_i,$ and the roots of $F'$ are all $\alpha_i.$ Now taking $F$ to a sufficiently high power proves the theorem.

The proof of Theorem 1 also bears some resemblance to the proof of the famous "Three Points Theorem" of G. Belyi \cite{Belyi}. In fact, similar ideas lead to a variation of its proof, obtained independently by G. Belyi himself, Leonardo Zapponi, and the author (see \cite{Belyi2}, \cite{Wolfart}, \cite{Zapponi}).
It is based on the following theorem.
\begin{Theorem} Suppose $n\ge 2$ and $q_1, q_2, ... ,q_n$ are rational numbers.
Then there exists a rational function with rational coefficients $f(x)$ with the following properties.
\begin{enumerate}
\item For every $i,$ $f(q_i)$ is either $0$ or $\infty$,
\item $f(\infty)=1,$
\item The map $f: P^1\rightarrow P^1$ is unramified outside $q_i$ and $\infty$.
\end{enumerate}
\end{Theorem}

{\bf Sketch of the proof.} We use the Lagrange interpolation formula for the constant polynomial $1$ at the points $\{q_i\}$. We multiply the resulting coefficients by a common denominator to get integers $\{k_i\}$. Then we choose $f(x)=\prod \limits_{i=1}^{n} (x-q_i)^{k_i}.$ 

{\bf Acknowledgments.} I became aware of the SR property and Thomas Craven's work on it from a question of Steve Fisk, posted on the sci.math.research newsgroup. I am greatly indebted to the people responsible for the creation and maintenance of this wonderful forum. I am also indebted to Bernhard K\"ock and Leonardo Zapponi for the references to \cite{Belyi2}, \cite{Wolfart} and \cite{Zapponi}, and to David Rohrlich for helpful suggestions.

\end{document}